\documentclass[12pt]{article}
\usepackage{amssymb,amsmath,amsthm,amsfonts}

\newcommand{\mS}{\overline{\mathcal{S}}}

\begin{document}

\begin{center}
 {\large \bf Reconstruction formula for differential systems with a singularity.
 } \\[0.2cm]
 {\bf Mikhail Ignatiev} \\[0.2cm]
\end{center}

\thispagestyle{empty}

{\bf Abstract.}
 Our studies concern some aspects of scattering theory of the singular differential systems
$
  y'-x^{-1}Ay-q(x)y=\rho By, \ x>0
$
with $n\times n$ matrices $A,B, q(x), x\in(0,\infty)$, where $A,B$ are constant and $\rho$ is a spectral parameter. We concentrate on the important special case when $q(\cdot)$ is smooth and $q(0)=0$ and derive a formula that express such $q(\cdot)$ in the form of some special contour integral, where the kernel can be written in terms of the Weyl - type solutions of the considered differential system.
Formulas of such a type play an important role in constructive solution of inverse scattering problems: use of such formulas, where the terms in their right-hand sides are previously found from the so-called main equation, provides a final step of the solution procedure.   In order to obtain the above-mentioned reconstruction formula we establish first the asymptotical expansions for the Weyl - type solutions as $\rho\to\infty$ with $o\left(\rho^{-1}\right)$ rate remainder estimate.
 \\
[0.1cm]

\section{Introduction}

Our studies concern some aspects of scattering theory of the differential systems
\begin{equation}\label{sys}
  y'-x^{-1}Ay-q(x)y=\rho By, \ x>0
\end{equation}
with $n\times n$ matrices $A,B, q(x), x\in(0,\infty)$, where $A,B$ are constant and $\rho$ is a spectral parameter.

Differential equations with coefficients having non-integrable singularities at the end or inside the interval often appear in various areas of natural sciences and engineering. For $n=2$, there exists an extensive literature devoted to different aspects of spectral theory of the radial Dirac operators, see, for instance {\cite{TKBDir}}, {\cite{AHM}}, {\cite{AHM1}}, {\cite{Ser}}, {\cite{GoYu}}.

Systems of the form {\eqref{sys}} with $n>2$ and arbitrary complex eigenvalues of the matrix $B$ appear to be considerably more difficult for investigation even in the "regular" case $A=0$ {\cite{BCsyst}}.  Some difficulties of principal matter also appear due to the presence of the singularity. Whereas the "regular" case $A=0$ has been studied fairly completely to date {\cite{BCsyst}}, {\cite{Zh}}, {\cite{Ysyst}}, for system {\eqref{sys}} with $A\neq 0$  there are no similar general results.

In this paper, we consider the important special case  when $q(\cdot)$ is smooth and $q(0)=0$ and, provided also that the discrete spectrum is empty, derive a formula that express such $q(\cdot)$ in the form of some special contour integral, where the kernel can be written in terms of the Weyl - type solutions of system \eqref{sys}.
Formulas of such a type play an important role in constructive solution of inverse scattering problems: use of such formulas, where the terms in their right-hand sides are previously found from the so-called \textit{main equation} (see, for instance, \cite{YIP93}, \cite{Ybook}), provides a final step of the solution procedure.   In order to obtain the above-mentioned reconstruction formula we establish first the asymptotical expansions for the Weyl - type solutions as $\rho\to\infty$ with $o\left(\rho^{-1}\right)$ rate remainder estimate.

\section{Preliminary remarks}

Consider first the following unperturbed system:
\begin{equation}\label{unp system}
  y'-x^{-1}Ay=\rho By
\end{equation}
and its particular case corresponding to the value $\rho=1$ of the spectral parameter
\begin{equation}\label{main system}
  y'-x^{-1}Ay=By
\end{equation}
but to {\it complex} (in general) values of $x$.

\textbf{Assumption 1.} Matrix $A$ is off-diagonal.
The eigenvalues $\{\mu_j\}_{j=1}^n$ of the matrix $A$ are distinct and such that $\mu_j-\mu_k \notin \mathbb{Z}$ for $j\neq k$, moreover, $\mbox{Re}\mu_1<\mbox{Re}\mu_2<\dots<\mbox{Re}\mu_n$, $\mbox{Re}\mu_k\neq 0$, $k=\overline{1,n}$.

\medskip
\textbf{Assumption 2.} $B=diag(b_1,\dots,b_n)$, the entries $b_1,\dots,b_n$ are nonzero distinct points on complex plane such that  $\sum\limits_{j=1}^n b_j=0$ and such that any 3 points are noncolinear.

\medskip
Under Assumption 1 system {\eqref{main system}}
has the fundamental matrix
$c(x)=(c_1(x),\dots,c_n(x))$, where
$$
c_k(x)=x^{\mu_k}\hat c_k(x),
$$
$\det c(x)\equiv 1$ and all $\hat c_k(\cdot)$ are entire functions, $\hat c_k(0)=\mathfrak{h}_k$, $\mathfrak{h}_k$ is an eigenvector of the matrix $A$ corresponding to the eigenvalue $\mu_k$. We define $C_k(x,\rho):=c_k(\rho x)$, $x\in(0,\infty)$, $\rho\in\mathbb{C}$. We note that the matrix $C(x,\rho)$ is a solution of unperturbed system {\eqref{unp system}} (with respect to $x$ for given spectral parameter $\rho$).

Let $\Sigma$ be the following union of lines through the origin in $\mathbb{C}$:
$$
\Sigma=\bigcup\limits_{(k,j): j\neq k}\left\{z:\mbox{Re}(z b_j)=\mbox{Re}(z b_k)\right\}.
$$
By virtue of Assumption 2 for any $z\in\mathbb{C}\setminus\Sigma$ there exists the ordering $R_1,\dots, R_n$ of the numbers $b_1,\dots,b_n$ such that $\mbox{Re}(R_1 z)<\mbox{Re}(R_2 z)\dots<\mbox{Re}(R_n z)$.
Let $\mathcal{S}$ be a sector $\{z=r\exp(i\gamma), r\in(0,\infty), \gamma\in(\gamma_1,\gamma_2)\}$ lying in $\mathbb{C}\setminus\Sigma$.
Then {\cite{Syb}} system \eqref{main system} has the fundamental matrix
 $e(x)=(e_1(x),\dots,e_n(x))$ which is analytic in $\mathcal{S}$, continuous in $\overline{\mathcal{S}}\setminus\{0\}$ and admits the asymptotics:
$$
e_k(x)=\mbox{e}^{ x R_k}(\mathfrak{f}_k+x^{-1}\eta_k(x)), \ \eta_k(x)=O(1),\ x\to\infty, \ x\in\overline{\mathcal{S}},
$$
where
$(\mathfrak{f}_1,\dots,\mathfrak{f}_n)=\mathfrak{f}$ is a permutation matrix such that $(R_1,\dots,R_n)=(b_1,\dots,b_n)\mathfrak{f}$. We define $E(x,\rho):=e(\rho x)$.

\medskip
Everywhere below we assume that the following additional condition is satisfied.

\medskip
\textbf{Condition 1.} For all $k=\overline{2,n}$ the numbers $$\Delta_{0k}:=\det(e_1(x),\dots,e_{k-1}(x),c_k(x),\dots,c_n(x))$$ are not equal to 0.

\medskip
Under Condition 1 system \eqref{main system} has the fundamental matrix $\psi_0(x)=(\psi_{01}(x),\dots,\psi_{0n}(x))$ which is analytic in $\mathcal{S}$, continuous in $\overline{\mathcal{S}}\setminus\{0\}$ and admits the asymptotics:
$$
\psi_{0k}(xt)=\exp(xt R_k)(\mathfrak{f}_k+o(1)), t\to\infty, x\in\mathcal{S} , \ \psi_{0k}(x)=O(x^{\mu_k}), x\to 0.
$$
We define $\Psi_0(x,\rho):=\psi_0(\rho x)$. As above, we note that the matrices $E(x,\rho)$, $\Psi_0(x,\rho)$ solve {\eqref{unp system}}.

\medskip
In the sequel we use the following notations:
\begin{itemize}
  \item $\{\mathfrak{e}_k\}_{k=1}^n$ is the standard basis in $\mathbb{C}^n$;
  \item $\mathcal{A}_{m}$ is the set of all ordered multi-indices $\alpha=(\alpha_1, \dots, \alpha_m)$, $\alpha_1<\alpha_2<\dots<\alpha_m$, $\alpha_j\in\{1,2,\dots,n\}$;
  \item for a sequence $\{u_j\}$ of vectors and a multi-index $\alpha=(\alpha_1,\dots,\alpha_m)$ we define $u_\alpha:=u_{\alpha_1}\wedge\dots\wedge u_{\alpha_m}$;
  \item for a numerical sequence $\{a_j\}$ and a multi-index $\alpha$ we define $$a_\alpha:=\sum\limits_{j\in\alpha} a_j, \ a^{\alpha}:=\prod\limits_{j\in\alpha}a_j;$$
  \item for a multi-index $\alpha$ the symbol $\alpha'$ denotes the ordered multi-index that complements $\alpha$ to $(1,2,\dots,n)$;
  \item for $k=\overline{1,n}$ we denote
$$\overrightarrow{a}_k:=\sum\limits_{j=1}^k a_j,\ \overleftarrow{a}_k:=\sum\limits_{j=k}^n a_j,\
\overrightarrow a^{k}:=\prod\limits_{j=1}^k a_j, \ \overleftarrow a^{k}:=\prod\limits_{j=k}^n a_j.$$
We note that Assumptions 1,2 imply, in particular,  $\sum\limits_{k=1}^n \mu_k=\sum\limits_{k=1}^n R_k=0$ and therefore for any  multi-index $\alpha$ one has $R_{\alpha'}=-R_\alpha$ and $\mu_{\alpha'}=-\mu_\alpha$.
  \item the symbol $V^{(m)}$, where $V$ is $n\times n$ matrix, denotes the operator acting in $\wedge^m \mathbb{C}^n$ so that for any vectors $u_1, \dots, u_m$ the following identity holds:
$$
V^{(m)}(u_1\wedge u_2 \wedge\dots\wedge u_m)=\sum\limits_{j=1}^m u_1\wedge u_2 \wedge\dots\wedge u_{j-1}\wedge V u_j\wedge u_{j+1}\wedge\dots\wedge u_m;
$$
  \item if $h\in \wedge^n \mathbb{C}^n$ then $|h|$ is a number such that $
h=|h|\mathfrak{e}_1\wedge\mathfrak{e}_2\wedge\dots\wedge\mathfrak{e}_n
$;
\item for $h\in \wedge^m \mathbb{C}^n$ we set:
$
\|h\|:=\sum\limits_{\alpha\in\mathcal{A}_m}|h_\alpha|,
$ where $\{h_\alpha\}$ are the coefficients from the expansion $h=\sum\limits_{\alpha\in\mathcal{A}_m}h_\alpha \mathfrak{e}_\alpha$.
\end{itemize}

\section{Asymptotics of the Weyl - type solutions}

Let $\mathcal S \subset \mathbb{C}\setminus\Sigma$ be an open sector with vertex at the origin. For arbitrary $\rho\in\mathcal S$ and $k\in\{1,\dots,n\}$ we define the \textit{$k$-th Weyl - type solution} $\Psi_k(x,\rho)$ as a solution of \eqref{sys} normalized with the asymptotic conditions:
\begin{equation}\label{psik def}
  \Psi_k(x,\rho) = O\left(x^{\mu_k}\right), \quad x\to 0, \qquad \Psi_k(x,\rho)=\exp(\rho x R_k)(\mathfrak f_k + o(1)), \quad x\to\infty.
\end{equation}
If $q(\cdot)$ is off-diagonal matrix function summable on the semi-axis $(0,\infty)$ then for arbitrary given $\rho\in\mathcal S$ $k$-th Weyl - type solution exists and is unique provided that the characteristic function:
$$
\Delta_k(\rho)=\left|F_{k-1}(x,\rho)\wedge T_k(x,\rho)\right|
$$
does not vanish at this $\rho$. Here $\left\{F_k(x,\rho)\right\}_{k=1}^n$, $\left\{T_k(x,\rho)\right\}_{k=1}^n$ are certain tensor-valued functions (\textit{fundamental tensors}) defined as solutions of certain Volterra integral equations, see \cite{Ign16}, \cite{IgnSSU19} for details.

For arbitrary fixed arguments $x,\rho$ (where $\Delta_k(\rho)\neq 0$) the value $\Psi_k=\Psi_k(x,\rho)$ is the unique solution of the following linear system:
\begin{equation}\label{psik ini sys}
  F_{k-1}\wedge \Psi_k = F_k, \quad \Psi_k\wedge T_k=0.
\end{equation}
This fact and also some properties of the Weyl - type solutions were established in works \cite{Ign16}, \cite{IgnMZ2020}, in particular, the following asymtotics for $\rho\to\infty$ was obtained:
\begin{equation}\label{psik as gen}
  \Psi_k(x,\rho) = \Psi_{0k}(x,\rho) + o\left(\exp(\rho x R_k)\right).
\end{equation}
For our purposes we need more detailed asymptocis that can be obtained provided that the potential $q(\cdot)$ is smooth enough and vanishes as $x\to 0$.

We denote by $\mathcal P(\mathcal S)$ the set of functions $F(\rho)$, $\rho\in\mathcal S$ admitting the representation:
$$
F(\rho)=\sum\limits_{\lambda\in\Lambda} f(\lambda)\exp(\lambda\rho).
$$
Here the set $\Lambda$ (depending on $F(\cdot)\in\mathcal P(\mathcal S)$) is such that $\mbox{Re}(\lambda\rho)<0$ for all $\lambda\in\Lambda$, $\rho\in\mathcal S$. We note that the set of scalar functions belonging to $\mathcal P(\mathcal S)$ is an algebra with respect to pointwise multiplication.

\medskip
\textbf{Theorem 1.} \textit{ Suppose that $q(\cdot)$ is absolutely continuous off-diagonal matrix function such that $q(0)=0$. Denote by $\hat q_o(\cdot)$ the off-diagonal matrix function such that $[B,\hat q_o(x)]=-q(x)$ for all $x>0$ (here $[\cdot,\cdot]$ denotes the matrix commutator). Define the diagonal matrix $d(x)=diag(d_1(x),\dots,d_n(x))$, where
$$d_k(x):=\int\limits_{x}^\infty t^{-1}\left(\left[\hat q_o(t),A\right]\right)_{kk} \,dt$$
and set $\hat q(x):=\hat q_o(x)+d(x)$.}

\textit{
Suppose that all the functions $q_{ij}(\cdot), q'_{ij}(\cdot)$ and $\tilde q_{ij}(\cdot)$, where $\tilde q(x):=\hat q'(x)+ x^{-1}[\hat q(x), A]$, belong to $X_p:=L_1(0,\infty)\cap L_p(0,\infty)$, $p>2$.}

\textit{
Then for each fixed $x>0$ and $\rho\to\infty$, $\rho\in\mS$ the following asymptotics holds:
$$
\rho(\Psi(q,x,\rho)-\Psi_0(x,\rho))\exp(-\rho x R)=\mathfrak{f}\Gamma(x)+\hat q(x)\mathfrak{f}+\mathcal{E}(x,\rho)+o(1),
$$
where $\Gamma(x)$ is some diagonal matrix, $\mathcal E(x,\cdot)\in\mathcal P(\mathcal S)$.}

\medskip
\textbf{Proof.} Denote:
$$
\tilde F_k(x,\rho):=\exp\left(-\rho x \overrightarrow R_k\right) F_k(x,\rho), \quad
\tilde T_k(x,\rho):=\exp\left(-\rho x \overleftarrow R_k\right) T_k(x,\rho).
$$
By virtue of Theorem 1 \cite{IgnSSU19}  the following asymptotics hold:
\begin{multline}\label{t as}
   \rho\tilde F_k(q,x,\rho)=\rho\tilde F_{0k}(x,\rho)+
   \sum\limits_{\alpha\in\mathcal A_k}f_{k,\alpha}(x)\mathfrak{f}_{\alpha}+\mathcal{E}(x,\rho)+o(1),\\
  \rho\tilde T_k(q,x,\rho)=\rho\tilde T_{0k}(x,\rho)+d_{0k}\tilde T_{0k}(x,\rho)+
  \sum\limits_{\alpha\in\mathcal A_{n-k+1}}T^0_{k,\alpha^*(k)}g_{k,\alpha,\alpha^*(k)}(x)\mathfrak{f}_{\alpha}+\mathcal{E}(x,\rho)+o(1),
\end{multline}
where $\alpha^*(k):=(k,\dots,n)$, $\alpha_*(k):=(1,\dots,k)$; $f_{k,\alpha}(x)$, $g_{k,\alpha,\alpha^*(k)}(x)$ are some scalars that can be written explicitly in terms of $q(\cdot)$.

For the Weyl - type solutions of the unperturbed system we have the asymptotics (following directly from their definition):
\begin{equation}\label{tilde psi}
  \tilde\Psi_{0k}(x,\rho)=\mathfrak{f}_k+\mathcal{E}(x,\rho)+O\left(\rho^{-1}\right),
\end{equation}
where $\tilde\Psi_{0k}(x,\rho):=\exp(-\rho x R_k)\Psi_{0k}(x,\rho)$.
Here and below we use the same symbol $\mathcal E(\cdot,\cdot)$ for different functions such that $\mathcal E(x,\cdot)\in\mathcal P(\mathcal S)$ for each fixed $x$.

We rewrite relations \eqref{psik ini sys} in the form of the following linear system with respect to value $\tilde\Psi_k=\tilde\Psi_k(x,\rho)$ of the function $\tilde\Psi_k(x,\rho):=\exp(-\rho x R_k)\Psi_k(x,\rho)$:
$$
\tilde F_{k-1}\wedge\tilde\Psi_k=\tilde F_k, \ \tilde\Psi_k\wedge\tilde T_k=0.
$$
By making the substitution:
\begin{equation}\label{anz}
  \tilde\Psi_k=\tilde\Psi_{0k}+\hat\Psi_k,
\end{equation}
we obtain:
$$
\tilde F_{k-1}\wedge\hat\Psi_k=\tilde F_k-\tilde F_{k-1}\wedge\tilde\Psi_{0k}, \
\hat\Psi_k\wedge\tilde T_k=-\tilde\Psi_{0k}\wedge \tilde T_k.
$$
The obtained relations we transform into the following system of linear algebraic equations:
\begin{equation}\label{SLAE}
  \sum\limits_{j=1}^n m_{ij}\gamma_{jk}=u_i, \ i=\overline{1,n}
\end{equation}
with respect to coefficients $\{\gamma_{jk}\}$ of the expansion:
\begin{equation}\label{exp Ga}
  \hat\Psi_k(x,\rho)=\sum\limits_{j=1}^n \gamma_{jk}(x,\rho) \mathfrak{f}_j.
\end{equation}
Coefficients $\{m_{ij}\}$, $\{u_i\}$ can be calculated as follows:
$$
m_{ij}=\left|\tilde F_{k-1}\wedge\mathfrak{f}_j\wedge\mathfrak{f}_\alpha\right|, $$$$ u_{i}=\left|(\tilde F_k-\tilde F_{k-1}\wedge\tilde\Psi_{0k})\wedge\mathfrak{f}_\alpha\right|, \ \alpha=\alpha^*(k)\setminus i, \quad i=\overline{k,n},
$$
$$
m_{ij}=\left|\mathfrak{f}_\alpha\wedge\mathfrak{f}_j\wedge\tilde T_k\right|, \
u_i=-\left|\mathfrak{f}_\alpha\wedge\tilde\Psi_{0k}\wedge\tilde T_k\right|, \ \alpha=\alpha_*(k-1)\setminus i, \quad i=\overline{1,k-1}.
$$
Using {\eqref{t as}}, {\eqref{tilde psi}} and taking into account that:
$$
\tilde F^0_{k-1}\wedge\tilde\Psi_{0k}=\tilde F^0_k, \ \tilde\Psi_{0k}\wedge\tilde T^0_k=0,
$$
we obtain the following asymptotics for the coefficients of SLAE \eqref{SLAE} as $\rho\to\infty$:
\begin{equation}\label{coe as big i}
  m_{ij}(x,\rho)=O\left(\rho^{-1}\right), j\neq i, \ m_{ii}(x,\rho)=m_{ii}^0+O\left(\rho^{-1}\right), \ m_{ii}^0=(-1)^{k-i}|\mathfrak{f}|,  \ i=\overline{k,n},
\end{equation}
and
$$
m_{ij}(x,\rho)=m^0_{ij}+O\left(\rho^{-1}\right), \ i=\overline{1,k-1},\ j=\overline{1,k-1},
$$
$$
m_{ij}(x,\rho)=m^0_{ij}+\mathcal E(x,\rho)+O\left(\rho^{-1}\right), \ i=\overline{1,k-1},\ j=\overline{k,n},
$$
where:
$$
m^0_{ij}=T^0_{k,\alpha^*(k)}|\mathfrak{f}_\alpha\wedge\mathfrak{f}_j\wedge\mathfrak{f}_{\alpha^*(k)}|,\ \alpha=\alpha_*(k-1)\setminus i,
$$
and therefore:
\begin{gather}\label{coe as small i}
  m_{ij}(x,\rho)=O\left(\rho^{-1}\right), \ j\neq i,\ j<k, \quad m_{ij}(x,\rho)=\mathcal E(x,\rho)+O\left(\rho^{-1}\right), \ j=\overline{k,n},\\
  m_{ii}(x,\rho)=m^0_{ii}+O\left(\rho^{-1}\right), \ m^0_{ii}=(-1)^{k-1-i}|\mathfrak{f}|T^0_{k,\alpha^*(k)}
  \end{gather}
for $i=\overline{1,k-1}$.

Proceeding in a similar way we obtain:
\begin{equation}\label{v}
  \rho u_i(x,\rho)=u^1_i(x)+\mathcal E(x,\rho)+o(1),
\end{equation}
\begin{equation}\label{v as big i}
u^1_i(x)=(-1)^{k-i}|\mathfrak{f}|f_{k,\alpha}(x)-\delta_{i,k}|\mathfrak{f}|f_{(k-1),\alpha_*(k-1)}(x), \quad
\alpha=\alpha_*(k-1)\cup\{ i\}, \quad i=\overline{k,n},
\end{equation}
where $\delta_{i,k}$ is a Kroeneker delta,
\begin{equation}\label{v as small i}
  u^1_i(x)=-(-1)^{k-i}\left|\mathfrak{f}\right|T^0_{k,\alpha^*(k)}g_{k,\beta,\alpha^*(k)}(x), \quad \beta=\alpha'\setminus k, \alpha=\alpha_*(k-1)\setminus i,
  \quad i=\overline{1,k-1}.
\end{equation}
Using the obtained asymptotics we obtain from {\eqref{SLAE}} the auxiliary estimate $\gamma_{ik}(x,\rho)=O\left(\rho^{-1}\right)$.

Then, using in {\eqref{SLAE}} the substitution $\gamma_{ik}(x,\rho)=\rho^{-1}\hat\gamma_{ik}(x,\rho)$ (where, as it was shown above, $\hat\gamma_{ik}(x,\rho)=O(1)$) we obtain for $i=\overline{k,n}$:
$$
m_{ii}(x,\rho)\hat\gamma_{ik}(x,\rho)=u_i^1(x)+\mathcal E(x,\rho)-\sum\limits_{j\neq i}m_{ij}(x,\rho)\hat\gamma_{jk}(x,\rho)+o(1).
$$
In view of {\eqref{coe as big i}}, {\eqref{v}} this yields:
\begin{equation}\label{hat gamma big i}
  \hat\gamma_{ik}(x,\rho)=\gamma^1_{ik}(x)+\mathcal E(x,\rho)+o(1), \ \gamma^1_{ik}=\frac{u^1_i(x)}{m^0_{ii}},
\end{equation}
$i=\overline{k,n}$.

Similarly, for $i<k$ we have:
$$
m_{ii}(x,\rho)\hat\gamma_{ik}(x,\rho)=u_i^1(x)+\mathcal E(x,\rho)-\sum\limits_{j\geq k}m_{ij}(x,\rho)\hat\gamma_{jk}(x,\rho)-
\sum\limits_{j<k, j\neq i}m_{ij}(x,\rho)\hat\gamma_{jk}(x,\rho)+o(1).
$$
Using {\eqref{coe as small i}} the obtained relation can be transformed as follows:
$$
m^0_{ii}\hat\gamma_{ik}(x,\rho)=u_i^1(x)+\mathcal E(x,\rho)-\sum\limits_{j\geq k}m_{ij}(x,\rho)\hat\gamma_{jk}(x,\rho)
+o(1).
$$
Now, using in the right hand side of the obtained formula {\eqref{coe as small i}} for $m_{ij}(x,\rho)$ and {\eqref{hat gamma big i}} for $\hat\gamma_{jk}(x,\rho)$ with $j=\overline{k,n}$ we conclude that formulas  {\eqref{hat gamma big i}} are true for $i<k$ as well.

In our further calculations we use particular form of the coefficients $f_{k,\alpha}(x)$ and $g_{k,\alpha,\beta}(x)$ given by Theorem 1 \cite{IgnSSU19}.

For $i=\overline{k,n}$ from {\eqref{hat gamma big i}}, {\eqref{v as big i}}, {\eqref{coe as big i}} we get:
\begin{equation}\label{gamma 1 big i}
  \gamma^1_{ik}(x)=\delta_{i,k}\tilde\gamma^1_{ik}(x)+f_{k,\alpha}(x), \ \alpha=\alpha_*(k-1)\cup i.
\end{equation}
Theorem 1 \cite{IgnSSU19} yields:
$$
f_{k,\alpha}(x)=\chi_\alpha\left|\left(\hat q^{(k)}(x)\mathfrak{f}_{\alpha_*(k)}\right)\wedge\mathfrak{f}_{\alpha'}\right|, \quad
\chi_\alpha:=|\mathfrak f_\alpha\wedge\mathfrak f_{\alpha'}|.
$$

Recall that any arbitrary linear operator $V$ acting in $\mathbb{C}^n$ can be expanded onto the wedge algebra $\wedge \mathbb{C}^n$ so that the identity
$$
V(h_1\wedge\dots\wedge h_m)=(Vh_1)\wedge\dots\wedge (Vh_m)
$$
remains true for any set of vectors $h_1,\dots, h_m$, $m\leq n$; moreover, for any $h\in\wedge^n\mathbb{C}^n$ one has $Vh=|V|h$  (here $|V|$ denotes determinant of matrix of the operator $V$ in the standard coordinate basis $\{\mathfrak{e}_1,\dots,\mathfrak{e}_n\}$). In what follows the symbol $\mathfrak{f}$ denotes the above mentioned expansion of the operator corresponding to the transmutation matrix $\mathfrak{f}$. We should note also that the relation $(\mathfrak{f}^{-1}V\mathfrak{f})^{(k)}=\mathfrak{f}^{-1}V^{(k)}\mathfrak{f}$ is true for any $n\times n$ matrix $V$. Taking this into account we obtain:
$$
f_{k,\alpha}(x)=\chi_\alpha\left|\left(\mathfrak{f}\left(\mathfrak{f}^{-1}\hat q^{(k)}(x)\mathfrak{f}\mathfrak{e}_{\alpha_*(k)}\right)\right)
\wedge\left(\mathfrak{f}\mathfrak{e}_{\alpha'}\right)\right|=
$$
$$
=|\mathfrak{f}_\alpha\wedge\mathfrak{f}_{\alpha'}||\mathfrak{f}|
\left|\left(\mathfrak{f}^{-1}\hat q^{(k)}(x)\mathfrak{f}\mathfrak{e}_{\alpha_*(k)}\right)
\wedge\mathfrak{e}_{\alpha'}\right|=
|\mathfrak{e}_\alpha\wedge\mathfrak{e}_{\alpha'}|
\left|\left(\left(\mathfrak{f}^{-1}\hat q(x)\mathfrak{f}\right)^{(k)}\mathfrak{e}_{\alpha_*(k)}\right)
\wedge\mathfrak{e}_{\alpha'}\right|.
$$
For the particular multi-index $\alpha=\alpha_*(k-1)\cup i$ arising at \eqref{gamma 1 big i} and arbitrary $n\times n$ matrix $V$ we have:
$$
|\mathfrak{e}_\alpha\wedge\mathfrak{e}_{\alpha'}|
\left|\left(V^{(k)}\mathfrak{e}_{\alpha_*(k)}\right)
\wedge\mathfrak{e}_{\alpha'}\right|=V_{ik}.
$$
Substituting the obtained relations into {\eqref{gamma 1 big i}} we arrive at:
\begin{equation}\label{gamma 1 big i fin}
  \gamma^1_{ik}(x)=\delta_{i,k}\tilde\gamma^1_{ik}(x)+\left(\mathfrak{f}^{-1}\hat q(x)\mathfrak{f}\right)_{ik}, \quad i=\overline{k,n}.
\end{equation}
Proceeding in a similar way in the case $i<k$, using {\eqref{coe as small i}}, {\eqref{v as small i}} we obtain:
\begin{equation}\label{gamma 1 small i}
  \gamma^1_{ik}(x)=g_{k,\beta,\alpha^*(k)}(x), \quad \beta=\alpha'\setminus k, \alpha=\alpha_*(k-1)\setminus i.
\end{equation}
Theorem 1 \cite{IgnSSU19} yields:
$$
g_{k,\alpha,\beta}(x)=\chi_\alpha\left|\left(\hat q^{(n-k+1)}(x)\mathfrak{f}_{\beta}\right)\wedge\mathfrak{f}_{\alpha'}\right|
$$
for $\beta\neq\alpha$. Repeating the arguments above we obtain:
$$
g_{k,\alpha,\beta}(x)=
|\mathfrak{e}_\alpha\wedge\mathfrak{e}_{\alpha'}|
\left|\left(\left(\mathfrak{f}^{-1}\hat q(x)\mathfrak{f}\right)^{(n-k+1)}\mathfrak{e}_{\beta}\right)
\wedge\mathfrak{e}_{\alpha'}\right|.
$$
In particular, one gets:
$$
g_{k,\beta,\alpha^*(k)}=
|\mathfrak{e}_\beta\wedge\mathfrak{e}_{\beta'}|
\left|\left(\left(\mathfrak{f}^{-1}\hat q(x)\mathfrak{f}\right)^{(n-k+1)}\mathfrak{e}_{\alpha^*(k)}\right)
\wedge\mathfrak{e}_{\beta'}\right|.
$$
If $\beta=\alpha'\setminus k$, $\alpha=\alpha_*(k-1)\setminus i$, $i<k$, then for arbitrary $n\times n$ matrix $V$ we have:
$$
|\mathfrak{e}_\beta\wedge\mathfrak{e}_{\beta'}|
\left|\left(V^{(n-k+1)}\mathfrak{e}_{\alpha^*(k)}\right)
\wedge\mathfrak{e}_{\beta'}\right|=V_{ik}.
$$
Substituting the obtained relations into {\eqref{gamma 1 small i}} we arrive at:
\begin{equation}\label{gamma 1 small i fin}
  \gamma^1_{ik}(x)=\left(\mathfrak{f}^{-1}\hat q(x)\mathfrak{f}\right)_{ik}, \quad i=\overline{1,k-1}.
\end{equation}
From {\eqref{gamma 1 small i fin}}, {\eqref{gamma 1 big i fin}}, {\eqref{hat gamma big i}} we obtain:
$$
\rho\gamma_{ik}(x,\rho)=\hat\gamma_{ik}(x,\rho)=\delta_{i,k}\tilde\gamma^1_{ik}(x)+\left(\mathfrak{f}^{-1}\hat q(x)\mathfrak{f}\right)_{ik}+\mathcal E(x,\rho)+o(1).
$$
In terms of the matrix $\gamma=(\gamma_{ik})_{i,k=\overline{1,n}}$ this is equivalent to:
$$
\rho\gamma(x,\rho)=\Gamma(x)+\mathfrak{f}^{-1}\hat q(x)\mathfrak{f}+\mathcal E(x,\rho)+o(1),
$$
where the matrix $\Gamma(x)$ is diagonal.
Finally, using {\eqref{exp Ga}} in the form $\hat\Psi(x,\rho)=\mathfrak{f}\gamma(x,\rho)$ we obtain the required relation.

$\hfil\Box$

\section{Reconstruction formula}

Let $\mathcal S_\nu$, $\nu=\overline{1,N}$ be the open pairwise nonintersecting sectors such that $\mathbb{C}\setminus\Sigma = \bigcup\limits_{\nu=1}^N \mathcal S_\nu$. Suppose that the sectors are enumerated in counterclockwise order. We denote by $\Sigma_\nu$ the open ray dividing $\mathcal S_\nu$ and $\mathcal S_{\nu+1}$ (assuming $\mathcal S_{N+1}:=\mathcal S_1$).
We agree that the rays $\Sigma_\nu$ are oriented $0$ to $\infty$. Denote by $\Sigma^+_\nu$ and $\Sigma^-_\nu$ the banks of cut (along $\Sigma_\nu$) belonging to $\mathcal S_{\nu+1}$ and $\mathcal S_\nu$ respectively. We agree that $\Sigma^+_\nu$ is oriented from $0$ to $\infty$ while
$\Sigma^-_\nu$ is oriented from  $\infty$ to $0$.

For a function $f(\rho)$, $\rho\in \mathcal S_\nu \cup \mathcal S_{\nu+1}$ and arbitrary $\rho_0\in\Sigma_\nu$ we denote by $f^\pm(\rho_0)$ the limit values (if they exist):
$$
f^-(\rho_0):= \lim\limits_{\rho\to\rho_0, \rho\in\mathcal S_\nu}f(\rho), \quad
f^+(\rho_0):= \lim\limits_{\rho\to\rho_0, \rho\in\mathcal S_{\nu+1}}f(\rho).
$$

We say that off-diagonal matrix function $q(\cdot)\in X_p$ belongs to the class $G_0^p$ if for any $\nu\in\{1,\dots,N\}$ and $k\in\{1,\dots,n\}$ it is true that $\Delta_k(\rho)\neq 0$ for all $\rho\in\overline{\mathcal S}_\nu$. If $q(\cdot)\in G^p_0$ then the limit values $\Psi^\pm_k(x,\rho_0)$ exist for any $k\in\{1,\dots,n\}$, $\rho_0\in\Sigma_\nu$, $\nu\in\{1,\dots,N\}$.

We denote by $\Psi(x,\rho)$  the matrix function $\Psi(x,\rho)=(\Psi_1(x,\rho),\dots,\Psi_n(x,\rho))$ and introduce the following \textit{spectral mappings matrix}:
$$
P(x,\rho):=\Psi(x,\rho)\Psi^{-1}_0(x,\rho).
$$
If $q(\cdot)\in G^p_0$ then the limit values $P^\pm_k(x,\rho_0)$ exist for any $k\in\{1,\dots,n\}$, $\rho_0\in\Sigma_\nu$, $\nu\in\{1,\dots,N\}$. Denote $\hat P(x,\rho):=P^+(x,\rho)-P^-(x,\rho)$. Following theorem contains main result of the paper.

\medskip
\textbf{Theorem 2.} \textit{Suppose that the potential $q(\cdot)\in G^p_0$ satisfies the conditions of Theorem 1. Then the following relation (reconstruction formula) holds:
$$
q(x)=\frac{1}{2\pi i}\int\limits_\Sigma \left[B,\hat P(x,\rho)\right] d\rho,
$$
where (as above) the brackets $[\cdot,\cdot]$ denote the matrix commutator and the integral is considered as the following limit (existing for each $x>0$):
$$
\frac{1}{2\pi i}\int\limits_\Sigma \left[B,\hat P(x,\rho)\right] d\rho:=\lim\limits_{r\to\infty}
\frac{1}{2\pi i}\int\limits_{\Sigma^r} \left[B,\hat P(x,\rho)\right] d\rho,
$$
$\Sigma^r:=\Sigma\cap\{\rho:|\rho|\leq r\}$.
}

\medskip
\textbf{Proof.}
Consider the function:
$$F(x,\rho):=\rho [B,P(x,\rho)]+q(x).$$

From Theorem 1 we have the asymptotics:
$$
\hat\Psi(x,\rho):=(\Psi(x,\rho)-\tilde\Psi(x,\rho))\exp(-\rho xR)=\rho^{-1}(\mathfrak{f}\Gamma_\nu(x)+\hat q(x)\mathfrak{f}+\mathcal E_\nu(x,\rho)+o(1))
$$
as $\rho\to\infty$, $\rho\in\mathcal S_\nu$,
where $R=diag(R_1,\dots,R_n)$, $\Gamma_\nu(x)$ are some diagonal matrices and $\mathcal E_\nu(x,\cdot)\in\mathcal P(\mathcal S_\nu)$.

For $\tilde\Psi_0(x,\rho)$ we have :
$$
\tilde \Psi_0(x,\rho)=\mathfrak{f}+\mathcal E_\nu(x,\rho)+o(1)
$$
as $\rho\to\infty$, $\rho\in\mathcal S_\nu$
(we use the same symbol for denoting possibly different functions from $\mathcal P(\mathcal S_\nu)$).

Since $|\det\tilde\Psi_0|=1$ the following asymptotics is also valid:
$$
\tilde \Psi^{-1}_0(x,\rho)=\mathfrak{f}^{-1}+\mathcal E_\nu(x,\rho)+o(1), \qquad\rho\to\infty, \quad \rho\in\mathcal S_\nu.
$$
Therefore, for $\rho\to\infty$, $\rho\in\mathcal S_\nu$ we have:
\begin{equation}\label{P as}
P(x,\rho)=I+\hat\Psi(x,\rho)\tilde\Psi^{-1}_0(x,\rho) =
I+\rho^{-1}(\mathfrak{f}\Gamma_\nu(x)\mathfrak{f}^{-1}+\hat q(x)+\mathcal E_\nu(x,\rho)+o(1)).
\end{equation}
Since the matrices $\Gamma_\nu(x)$ are diagonal the matrices $\mathfrak f\Gamma_\nu(x)\mathfrak f^{-1}$ are diagonal as well and we have $[B,\mathfrak f\Gamma_\nu(x)\mathfrak f^{-1}]=0$. Thus, from \eqref{P as} we deduce:
\begin{equation}\label{F ass}
  F(x,\rho)=\mathcal E_\nu(x,\rho)+o(1), \qquad\rho\to\infty, \quad \rho\in\mathcal S_\nu.
\end{equation}
Define:
$$
\gamma = \bigcup\limits_{\nu=1}^N \left(\Sigma^-_\nu\cup\Sigma^+_\nu\right), \qquad
\gamma_r:=\gamma\cap\{\rho:|\rho|\leq r\}, \quad \Gamma_r:=\gamma_r\cup C_r,
$$
where $C_r$ is the circle $\{\rho:|\rho|=r\}$ with a counterclockwise orientation.

By virtue of the Jordan lemma from asymptotics {\eqref{F ass}} it follows that for any arbitrary fixed $\rho\in\mathbb{C}\setminus\Sigma$ we have:
$$
\lim\limits_{r\to\infty}\int\limits_{C_r} \frac{d\zeta}{\zeta-\rho} F(x,\zeta)=0.
$$
Therefore, the Cauchy integral formula for the closed contour $\Gamma_r$ (where $r>|\rho|$):
$$
F(x,\rho)=\frac{1}{2\pi i}\int\limits_{\Gamma_r} \frac{d\zeta}{\zeta-\rho}F(x,\zeta)
$$
can be transformed as follows:
$$
F(x,\rho)=\lim\limits_{r\to\infty}\frac{1}{2\pi i}\int\limits_{\Sigma^r} \frac{d\zeta}{\zeta-\rho}(F^+(x,\zeta)-F^-(x,\zeta)).
$$
Taking into account that $F^+(x,\zeta)-F^-(x,\zeta)=\zeta[B,\hat P(x,\zeta)]$ we obtain:
\begin{equation}\label{F Cauchy}
  F(x,\rho)=\lim\limits_{r\to\infty}\frac{1}{2\pi i}\int\limits_{\Sigma^r} \frac{d\zeta}{\zeta-\rho}\zeta[B,\hat P(x,\zeta)].
\end{equation}
On the other hand, we can proceed in a similar way applying the Cauchy formula to the function $P(x,\rho)-I$. Thus we obtain:
$$
P(x,\rho)-I=\frac{1}{2\pi i}\int\limits_{\Gamma_r} \frac{d\zeta}{\zeta-\rho}(P(x,\zeta)-I)
$$
and since from {\eqref{F ass}} it follows that:
$$
\lim\limits_{r\to\infty}\int\limits_{C_r} \frac{d\zeta}{\zeta-\rho} (P(x,\zeta)-I)=0
$$
we arrive at the representation:
$$
P(x,\rho)=I+\lim\limits_{r\to\infty}\frac{1}{2\pi i}\int\limits_{\Sigma^r} \frac{d\zeta}{\zeta-\rho}(P^+(x,\zeta)-P^-(x,\zeta)).
$$
Substituting this to the definition of the function $F(x,\rho)$ we arrive at the representation:
$$
F(x,\rho)=q(x)+\lim\limits_{r\to\infty}\frac{1}{2\pi i}\int\limits_{\Sigma^r} \frac{d\zeta}{\zeta-\rho}\rho[B,\hat P(x,\zeta)].
$$
Compare it with {\eqref{F Cauchy}} we obtain the desired relation.

$\hfil\Box$

\medskip
{\bf Acknowledgements.} This work was supported by the Russian Science Foundation (project no. 19-01-00102, 20-31-70005).

\noindent Ignatiev, Mikhail\\
Department of Mathematics, Saratov State University, \\
Astrakhanskaya 83, Saratov 410012, Russia, \\
e-mail: mikkieram@gmail.com, ignatievmu@sgu.ru

\end{document}